\documentclass[12pt]{article}
\usepackage{graphicx}
\begin{document}
\vskip+0.2cm

\begin{center}
	{\bf\Large The Numerical Solution of the External Dirichlet Generalized Harmonic Problem for a Sphere by the Method of Probabilistic Solution}\\[5mm]
	
	{\sc \bf Mamuli Zakradze$^1$, Zaza Tabagari$^1$, Nana Koblishvili$^1$, Tinatin Davitashvili$^2$, Francisco Criado-Aldeanueva$^3$}\\[3mm]
	
	{$^1$Department of Computational Methods, Muskhelishvili Institute of Computational Mathematics, Georgian Technical University}\\[2mm]
	{email: mamuliz@yahoo.com , z.tabagari@hotmail.com, nanakoblishvili@yahoo.com}\\[2mm]
	
	{$^2$Iv.Javakhishvili Tbilisi State University, FENS}\\[2mm]
	{email: tinatin.davitashvili@tsu.ge}\\[2mm]

	{$^3$Department of Applied Physics, II, Polytechnic School, Malaga University}\\[2mm]
	{email: fcaldeanueva@ctima.uma.es}
\end{center}

{\bf{Abstract}}\\
\


 In the present paper an algorithm for the numerical solution of  the external
 Dirichlet generalized harmonic problem for a sphere by the method of probabilistic
 solution (MPS) is given. Under a generalized problem is meant the problem when a boundary function  has a finite number of curves with the first kind of discontinuity. The algorithm consists of the following main stages: 1) transition from an infinite domain to a finite domain by an inversion; 2) consideration of a new Dirichlet generalized harmonic problem on the basis of Kelvin's theorem for the obtained finite domain; 3) application of the MPS for the numerical solution a new problem which in turn is based on a computer simulation of the Wiener process; 4) finding the probabilistic solution of the posed generalized problem at any fixed points of the infinite domain by the solution of the new problem. For illustration numerical examples are considered and results are presented.

\\
 \\
 {\bf{1. Introduction}}\\
 \

First of all, we should note that the requirement of continuity for the boundary function in Dirichlet's classical harmonic problem is a very strong constraint, because in practical stationary problems (connected with electric, thermal and other static fields), there are cases when it is necessary to discuss and study $2D$ or $3D$ Dirichlet generalized harmonic problems. The problems of this type appeared in the literature mainly from the 40-s of the 20th century (see, e.g., [1-5]).

It is known (see e.g., [1,6]) that the methods used to obtain an approximate solution of the classical  boundary-value problems are: a) less suitable or b) not suitable at all for solving generalized boundary value problems. In the first case  convergence of corresponding approximate process is very slow in neighborhood of boundary singularities and, consequently, the accuracy of  approximate solution of the generalized problem is very low (see, e.g., [1-5]). In the second case, the process is unstable. For example, a similar phenomenon takes place when solving the 3D Dirichlet generalized harmonic problem by the method of fundamental solutions. 

In the above-mentioned literature, simplified, or so-called "solvable" generalized problems (the problems "whose" solutions can be constructed by series with terms, represented by special functions) are considered. The methods of separation of variables, particular solutions, and heuristic methods are mainly applied for their solution, therefore the accuracy of the solution is rather low. Since heuristic methods do not guarantee to find the best solution (moreover, in some cases they may give an incorrect solution), it is necessary to check these solutions in order to establish how well they satisfy all the conditions of a problem (see e.g., [1]). Besides, it should be noted that in the literature (see,e.g., [4], pp.346-348), while solving the Dirichlet external generalized harmonic simplest problem for a sphere, the existence of discontinuity curve is ignored.

The study of Dirichlet's generalized $2D$ and $3D$ problems in terms of the existence and uniqueness of solution, and selection of a reliable and effective method for its numerical solution has been intensively carried out since the 21st century in the Department of Computational Methods of the Niko Muskhelishvili Institute of Computational Mathematics.

The choice and construction of computational schemes (algorithms) mainly depend on the problem class, dimension, geometry and location of singularities on the boundary.
In particular, for solving the Dirichlet generalized plane harmonic problems the following approaches may be used: I) A method of reduction of Dirichlet generalized harmonic problems to a classical problem (see, e.g., [7, 8]); II) A method of conformal mapping (see, e.g., [9]); III) A method of probabilistic solution (see, e.g., [10, 11]). It is evident, that in the case of $3D$ Dirichlet harmonic problems, from the above mentioned approaches we can apply only third one.

In order to imagine the difficulties associated with the numerical solution of not only generalized, but even the classical external 3D harmonic Dirichlet problem, we will give below a brief overview of several works related to these topics.

In [12], the boundary conditions for the numerical solution of elliptic equations in exterior regions are considered. Such equations in exterior regions frequently require a boundary condition at infinity to ensure the well-posedness of the problem. Practical applications include examples of the Helmholtz and Laplace’s equations. The constructed algorithm requires the replacement of the condition at infinity by a boundary condition on a finite artificial surface. Direct imposition of the condition at infinity along the finite boundary results in large errors. A sequence of boundary conditions is developed which provides increased accuracy of approximations to the problem in the infinite domain. Estimates of the error is obtained for several cases.

In [13], the Laplace-Dirichlet problem is investigated in three-dimensional case. The Laplace-Dirichlet problem is read in the same words as the Laplace-Neumann equation, excepted for the boundary condition. Corresponding variational formulation is considered. It based on the introduction of an auxiliary unknown by the means of  decomposition of the function $u$.
 
In the paper [14], the analytic-numerical method for solving 3D exterior problems for elliptic equations under Dirichlet and Neumann boundary conditions in half-space is considered. Based on analytical transformation, the external boundary problem is reduced to the internal one, then the corresponding difference problem is considered based on the grid methods [6-8]. The implementation procedure of the considered analytical-numerical method for solving external boundary problems in three-dimensional semi-space reduces the problem to such a problem that it’s possible using traditional techniques and methods of numerical analysis.

In [15] is considered the solution of direct and inverse exterior boundary value problems via the strongly conditioned stochastic method, mainly for exterior harmonic problems. Note that in numerical examples from the Section 4.2, was performed experiments with cones of different thickness and observation points in the near and far field. But, actually the observation points given in Table 2 are “not too far” from the field, when we are considering exterior boundary value problems.

It is known that elliptic problems in external domains arise in many branches of physics. For example, the Laplace equation $(\Delta u(x)=0, x \in D )$ arises in the studies of thermostatic and electrostatic fields external to given surfaces(see, e. g.,[3,16]); the flow of an incompressible irrotational fluid around a body is described by the same equation (see, e. g.,[17]) and so on.

In these cases infinity can be regarded as a separate boundary. A condition at infinity is required to make the external problem well posed. For the Laplace equation it is sufficient to impose a condition of regularity at infinity. In three dimensional $(3D)$ case the condition is $\lim u(x)=0$, for $|x| \rightarrow \infty$, where $|x|$ is the distance from a fixed (but chosen arbitrary) origin (see, e.g.,[18,19]).

Let $D=R^{3}\backslash K(O,R)$ be an infinite domain in the Euclidean space $R^{3}$, where $K(O,R)\equiv K_{R}$ is a kernel with a spherical surface $S(O,R)\equiv S_{R}$ (with the center at the point $O$ and the radius $R$, respectively). Since the harmonicity of a function is invariable under the linear transformation of the Cartesian coordinate system, therefore without loss of generality we assume that the origin of coordinates is at the point $O=(0,0,0)$.

Let us formulate the following problem.

{\bf{Problem A.}} {\it{The function $g(y)$ is given on the boundary $S_{R}$ of the infinite region  $D$ and is continuous everywhere, except a finite number of curves $l_{1}, l_{2}, ...,l_{n}$, which represent discontinuity curves of the first kind for the function $g(y)$. It is required to find a function
 $u(x)\equiv u(x_{1},x_{2},x_{3})\in C^{2}(D)\bigcap C(\overline {D}\backslash\bigcup\limits_{k=1}^{n}
 l_{k})$ satisfying the conditions:
$$
\Delta u(x)=0,\;\;\;\;\;x\in D,  \eqno(1.1)
$$
$$
u(y)=g(y),\;\;\;\;\; y \in S_{R},\;\;\;y\; \overline {\in}\; l_{k}\subset S_{R} \;\;(k=\overline {1,n}), \eqno (1.2)
$$
$$
|u(x)|<c,\;\;\;\;\;\;\;\;\;\;\;\;\;x\in \overline {D},\;\eqno (1.3)
$$
$$
\lim u(x)=0\;\;for\;\;x\rightarrow \infty,\;\eqno (1.4)
$$
where $\Delta=\sum\limits_{i=1}^3\frac{\partial^2}{\partial x^2_i}$ is the Laplace
operator, $c$ is a real constant.}}

On the basis of (1.3), in general, the values of $u(y)$ are not uniquely defined on the curves $l_{k} (k=\overline {1,n})$. In particular, if Problem $A$ concerns the determination of the thermal (or the electric) field, then we must take $u(y)=0$ when $y\in l_{k}$, respectively. In this case, in the physical sense the curves $l_{k}$ are non-conductors (or dielectrics). Otherwise, $l_{k}$ will not be a discontinuity curve.

It is evident that, surface $S_{R}$ is divided into the parts $S_{R}^{i} (i=\overline {1,m})$ by the curves $l_{k} (k=\overline {1,n})$, where one of the following conditions holds: $n=m, n<m, n>m$. Thus, the boundary function $g(y)$ has the following form
$$
g(y)=\left \{\begin {array}{ll}
g_1(y),\;\;\;y\in S_{R}^1,\\&\\
g_2(y),\;\;\;y\in S_{R}^2,\\&\\
\ldots \ldots \ldots \ldots \\&\\
g_{m}(y),\;\;\;y\in S_{R}^{m},\\&\\
\ 0,\;\;\;y\in l_{k}\; (k=\overline{1,n}),
\end {array}\right .
\eqno (1.5)
$$
where: $ S_{R}^i\; (i=\overline{1,m})$ are the parts of $S_R$ without
discontinuity curves, respectively; the functions $g_{i}(y),\; y\in
S_{R}^i \; (i=\overline{1,m} )$ are continuous on the parts $S_{R}^i $.
It is evident that $ S_{R}=(\bigcup \limits _{j=1}^{m}S_{R}^{j})\bigcup
(\bigcup \limits_{k=1}^{n}l_{k}).$
\\
\\
{\bf{Remark 1.}} If the domain $D$ is finite with a surface $S$, then the problem of type $A$ with the conditions (1.1), (1.2), (1.3) has a unique solution depending continuously on the data (see [20,21]).
\\
\\
It should be noted that the condition (1.4) is essential for the uniqueness of solution of the Problem $A$. Indeed, if out of (1.4) Problem $A$ has any solution $u_{1}(x)$, then $u_{2}(x)=u_{1}(x)+k(r-R)/r$ is its solution as well, where $r\geq R$ and $k$ is a real constant, i.e., it has an infinite set of solutions.

On the other hand, if Problem $A$ has a solution, then it is unique. Assume that Problem $A$ with the boundary function (1.5) has two solutions $u_{1}(x)$ and $u_{2}(x)$, then $u(x)=u_{1}(x)-u_{2}(x)$ is a solution of Problem $A$ with null boundary value and for function $u(x)$ condition (1.4) is
fulfilled. It is known (see e.g.,[19], p.303) that, in the noted case $u(x)\equiv 0 $ when
$x\in \overline D$ (or $ u_{1}(x)=u_{2}(x), x\in \overline D)$. The existence of solution of the Problem $A$ will be shown in Section 2. Based on the above-mentioned, we can investigate by Problem $A $ only such physical phenomena, which are damped at infinity.

In the case of external 3D Dirichlet generalized harmonic problems, the difficulties become more significant. In particular, there does not exist a standard algorithm which can be applied to a wide class of domains.

On the basis of noted above, for the numerical solution of the Problem $A$ we should apply an algorithm which does not require the approximation of a boundary function and in which the existence of discontinuity curves is not ignored. The method of probabilistic solution (MPS) is one such method
(see [20,21]), but its direct application in infinite domains is impossible.

Performed investigations showed (see, e.g., [10,11,22-27]) that the MPS is ideally suited for numerical solution of both classical and generalized (2D and 3D) Dirichlet harmonic problems for a wide class of finite domains only.

Therefore, construction of efficient high accuracy computational schemes for the numerical solution of external 3D Dirichlet generalized harmonic problems (which can be applied to a wide class of domains) has both theoretical and practical importance.
\\
\\
{\bf {2. Transition from the infinite domain $D$ to the kernel $K_{R}$ by an inversion
and consideration of a new problem on the basis of Kelvin's theorem}}
\\

In the Problem $A$ the domain $D$ is infinite, therefore it is impossible the direct application of the MPS to its solving. In order to solve the Problem $A$ using MPS, we perform transition from the domain $D$ to the kernel $K_{R}$ by means of the inversion (see e.g., [4,28]). Let a point $x(x_1,x_2,x_3) \in D$ and consider the following inversions:
$$
\xi_{i}=\frac{R^{2}}{|x|^{2}}x_{i}, \; x\in D, \; (i=1,2,3),\; |x|^2=x_{1}^{2}+x_{2}^{2}+x_{3}^{2}, \eqno (2.1)
$$
$$
x_{i}=\frac{R^{2}}{|\xi|^{2}}\xi_{i}, \; \xi \in K_{R}\backslash \{O\},\; (i=1,2,3),\; |\xi|^2=\xi_{1}^{2}+\xi_{2}^{2}+\xi_{3}^{2}, \eqno (2.2)
$$
with respect to the sphere $S_{R}$. The points $x$ and $\xi$ are called symmetric points with respect to the sphere $S_{R}$. From (2.1) (or (2.2)) we have
$$
|x||\xi|=R^{2},\; |\xi| \neq 0.   \eqno (2.3)
$$

It is known (see e.g., [4], p.260) that on the basis of (2.1) (or (2.2)) the symmetric points $\xi$ and $x$ with respect to the sphere $S_{R}$ are situated on the ray, whose beginning is at the point $|\xi|=0$ $(or\; \xi\equiv O)$.
It is easy to see that by (2.3) the infinite domain $ D$ is transformed one-to-one onto the $K_{R}$. In particular, the points of $S_{R}$ are transformed into itself, and the point $x=\infty$ is transformed into the point $\xi=O$ and vice versa.

On the basis of (2.2) the functions $u(x)$ and $g(y)$ are transformed into the functions
$$
u(\xi)\equiv u(\frac{R^{2}}{|\xi|^{2}}\xi_{1},\frac{R^{2}}{|\xi|^{2}}\xi_{2},\frac{R^{2}}{|\xi|^{2}}\xi_{3})\;\; and \;\; g(\eta)\equiv g(\frac{R^{2}}{|\eta|^{2}}\eta_{1},\frac{R^{2}}{|\eta|^{2}}\eta_{2},\frac{R^{2}}{|\eta|^{2}}\eta_{3}), \eqno (2.4)
$$
respectively, where $\xi(\xi_{1},\xi_{2},\xi_{3})\in K_R \backslash S_{R},\; |\xi|\neq 0
, \;\;\eta(\eta_{1},\eta_{2},\eta_{3})\in S_{R}.$

It is easy to see that the function $ u(\xi)$ is not harmonic in the domain $ K_{R}\backslash S_{R}$. But, we can remove the noted defect, if we apply Kelvin's theorem [4,28].

{\bf{Theorem 1.}} {\it{If a function $u(x_{1},x_{2},x_{3})$ is harmonic in the infinite domain $D$, then the function
$$
v(\xi_{1},\xi_{2},\xi_{3})=\frac{R}{|\xi|}u(\xi) \eqno (2.5)
$$
is harmonic in the domain $K_{R}\backslash \{O\}$.}}

The function $v(\xi)$ is bounded in the neighborhood of the point $\xi^{0}=(0,0,0)$,
therefore this point for $v(\xi)$ is a removable singular point [4,28]. But, actually, on the basis of the Theorem 1, for extension of harmonicity of function $v(\xi)$ at the point $\xi^{0}$ we must solve the following Dirichlet generalized harmonic problem.

{\bf{Problem A*.}} {\it{Find a generalized harmonic function $v(\xi)$ satisfying the conditions:
$$
\Delta v(\xi)=0,   \xi\in K_{R}\backslash S_{R},   \eqno (2.6)
$$
$$
v(\eta)=g(\eta),\; \eta \in S_{R},\; \eta \overline{\in}\; l_{k}\subset S_{R},\;  \eqno (2.7)
$$
$$
|v(\xi)|<c,\;\;\xi \in K_{R},  \eqno (2.8)
$$
where in (2.7): $|\eta|=R, \ \eta \equiv y\  (y\in S_{R}), \ l_{k}\  (k=\overline{1,n})$ and $c$ are the same that in Problem $A$.}}

It is known (see [20,21]) that the generalized problem (2.6),\ (2.7),\ (2.8) is correct, i.e., the solution exists, is unique and depends continuously on the data. Respectively, Problem $A$ is correct.

It is evident that for solving the Problem $A^{*}$ we can apply the MPS. In particular, if we want to find the value of the solution $u(x)$ of Problem $A$ at a point $x (x\in D)$, first of all we have to find the image $\xi$ of $x$ by means of (2.1) and then find the solution $v(\xi)$ of the Problem $A^{*}$ at the point $\xi$.
Finally, on the basis of (2.5) we have
$$
u(x)=\frac{v(\xi)}{|x|}R,    \eqno (2.9)
$$
where $x\in D, \xi\in K_{R}\backslash S_{R}.$ 

Thus, for definition of value of the solution to Problem $A$ at the point $x$ we have
the formula (2.9). It is easy to see that for the function $u(x)$, defined by (2.9),
the conditions of Problem $A$ are fulfilled.
\\
\\
{\bf {3. The method of probabilistic solution and simulation of the Wiener process}}
\\

In this section we describe the solution of Problem $A^{*}$ by the MPS.
It is known(see e.g.,[21,24]), that the probabilistic solution of Problem $A^{*}$ at the fixed point $\xi \in K_{R}\backslash S_{R}$ has the following form
$$
v(\xi)=E_{\xi}g(\xi(\tau )). \eqno (3.1)
$$

In (3.1) $E_{\xi} g(\xi(\tau))$ is the mathematical expectation of values of the boundary function $g(\eta)$ at the random intersection
points of the trajectory of the Wiener process and the boundary $S_{R}$;
$\tau$ is a random moment of the first exit of the Wiener process
$\xi(t)=(\xi_1(t),\xi_2(t),\xi_3(t))$ from  the domain $K_{R}$. It is assumed
that the starting point of the Wiener process is always
$\xi(t_0)=(\xi_1(t_0),\xi_2(t_0),\xi_3(t_0))\in K_{R}\backslash S_{R}$ , where the value of desired function is being determined. If the number $N$ of the
random intersection points $\eta^j=(\eta_1^j,\eta_2^j,\eta_3^j)\in S_{R}\ \
(j=\overline{1,N})$ is sufficiently large, then according to the law
of large numbers, from (3.1) we have
$$
v(\xi) \approx  v_N(\xi)=\frac{1}{N}\sum\limits _{j=1}^N g(\eta^j) \ \
\eqno(3.2)
$$
or $v(\xi)=\lim v_N(\xi)$ for $N\rightarrow \infty$, in probability.
Thus, if we have the Wiener process, the approximate value of the
probabilistic solution to the Problem $A^{*}$ at point $\xi\in K_{R}\backslash S_{R}$ is calculated by the formula (3.2).

 In order to simulate the Wiener process we use the following
recursion relations (see e.g.,[20,24]):
$$
\xi_1(t_k)=\xi_1(t_{k-1})+\gamma_1(t_k)/nq,
$$
$$
\xi_2(t_k)=\xi_2(t_{k-1})+\gamma_2(t_k)/nq,\ \ \ \eqno(3.3)
$$
$$
\xi_3(t_k)=\xi_3(t_{k-1})+\gamma_3(t_k)/nq,
$$
$$
(k=1,2,\cdots),\ \xi(t_0)=\xi,
$$
according to which the coordinates of point $\xi(t_k) =
(\xi_1(t_k), \xi_2(t_k), \xi_3(t_k))$ are being determined. In (3.3)
$\gamma_1(t_k),\gamma_2(t_k),\gamma_3(t_k)$ are three normally
distributed independent random numbers for the $k$-th step, with
zero means and variances equal to one (The above numbers are generated apart);
 $nq$ is the number of quantification $(nq)$ such that
 $1/nq =  \sqrt{t_k-t_{k-1}}$ and when $nq\rightarrow\infty$, then the discrete
process approaches the continuous Wiener process. In the implementation, random process is simulated at each step of the walk and continues until its trajectory crosses the boundary.

In the considered case computations are performed and the random numbers are generated in MATLAB.

{{%
}}
{{%
}}
{\bf {Remark 2. In general, problems of type
{A}

can be solved by the MPS for all such locations of discontinuity curves, which give the possibility to establish the part of surface {S} where the intersection point is located.}}

{{%
}}
{{%
}}
{\bf {4. Numerical examples}}
{{%
}}
{{
}}It should be noted that in the 3D case there are no test solutions for
generalized problems of type $A^{*}$, therefore, for the verification of the
scheme needed for the numerical solution of Problem $A^{*}$, the reliability of
obtained results can be demonstrated in the following way.

If in boundary conditions (2.7) of Problem $A^{*}$ we take
$g_i(\eta)=1/|\eta-\xi^0|$, where $\eta\in S_{R}^i (i=\overline {1,m}),
\xi^0=(\xi_1^0,\xi_2^0,\xi_3^0)\ \overline \in\ K_{R}$, and $|\eta-\xi^0|$
denotes the distance between points $\eta$ and $x^0$, then it is
evident that the curves $l_k (k=\overline {1,n})$ represent
removable discontinuity curves for the boundary function $g(\eta)$.
Actually, in the mentioned case instead of generalized problem $A^{*}$
we obtain the following Dirichlet classical harmonic problem.

{\bf  {Problem B.}}{\it { Find a Function $v(\xi)\equiv v(\xi_1,\xi_2,\xi_3)
\in C^2(K_{R}\backslash S_{R} )\bigcap C(K_{R})$ satisfying the conditions:
$$
\Delta v(\xi)=0,\;\; \xi \in K_{R}\backslash S_{R},
$$
$$
v(\eta)=1/|\eta-\xi^0|,\;\; \eta \in S_{R},\; \xi^0 \;  \overline \in \; K_{R}.
$$}}
We solve this problem (by the MPS) using the program used for the Problem $A^{*}$. It is well-known that the Problem $B$ is well posed, i.e., its solution exists, is unique and depends on data continuously. Evidently, an exact solution of the Problem $B$ is
$$
v(\xi^0,\xi)=\frac{1}{|\xi-\xi^0|},\;\xi \; \in \; K_{R},\; \xi^0\;
\overline \in \;K_{R}. \eqno (4.1)
$$

Note that application of the MPS for numerical solution of the Dirichlet classical
harmonic problems is interesting and important (see
e.g.,[10, 11, 22, 23]). In this paper, the Problem $B$ has an auxiliary role.
In particular, for the Problem $B$, verification of the scheme needed
for the numerical solution of Problem $A^{*}$ and corresponding
program (comparison of the obtained results with exact solution) are
carried out first of all, and then actually Problem $A^{*}$ is solved under the
boundary conditions (1.5).

In the present paper MPS is applied for two examples. In the tables,
$N$ is a number of implementation of the Wiener process for the
given points $\xi^i=(\xi_1^i,\xi_2^i, \xi_3^i)\in K_{R}\backslash S_{R}$, and $nq$ is a number of quantification. For simplicity, in the considered examples
the values of $nq$ and $N$ are the same. In the tables for problems of
type $B$ the absolute errors $\Delta^i $ at the points
$\xi^i\in K_{R}\backslash S_{R}$ of $v_N(\xi)$ are presented in the MPS approximation, for $nq=100$ and various values of $N$. In particular, we have $\Delta^i=|v_N(\xi^i)-v(\xi^0,\xi^i)|$, where $v_N(\xi^i)$ is the
approximate solution of Problem $B$ at the point $\xi^i$, which is
defined by formula (3.2), and $v(\xi^0,\xi^i)$ is an exact solution of the
test problem is given by (4.1). In tables, for the problems of type $A^{*}$,
the probabilistic solution $v_N(\xi)$, defined by (3.2), is presented at the points
$\xi^i\in K_{R}\backslash S_{R}$.
\\
\\
{\bf {Remark 3.}} The problems of type $A^{*}$ and $B$ for ellipsoidal,
spherical, cylindrical, conic, prismatic, pyramidal and axisymmetric finite domains with a cylindrical hole are considered in [24-27].
\\
\\
{\bf{Example 1}}. In the first example exterior of the unit sphere $S_{1}:$ $y_{1}^{2}+y_{2}^{2}+y_{3}^{2}=1$, with the center at origin $O(0,0,0)$ and radius $R=1$ is considered in the role of domain $D$, where
$(y_{1},y_{2},y_{3})$ is a point of the surface $S_{1}$. In the considered case, in Problem $A$ the boundary function $g(y)$ has the following form
$$
g(y)=\left \{\begin {array}{ll}
1,\;\;\;y\in S^1=\{y\in S_{1}|\;-1\leq y_{3} <-0.5 \},\\&\\
2,\;\;\;y\in S^2=\{y\in S_{1}|\;-0.5<y_{3}<0.5\},\\&\\
1,\;\;\;y\in S^3=\{y\in S_{1}|\;0.5<y_{3}\leq 1\},\\&\\
0,\;\;\;y\in l_k\  (k=\overline {1,2}).\\&\\
\end {array}\right .
\eqno (4.2)
$$

 It is evident that in (4.2), the discontinuity curves $l_1,l_2$ are the circles, which are obtained by intersection of the planes $x_{3}=-0.5, x_{3}=0.5$,
 and the surface $S_{1}$ (in the physical sense curves $l_{1}$ and $l_{2}$ are non-conductors).

 According to above noted, for solving the Problem $A$, in the first place, we solve Problem $A^{*}$  with boundary function $g(\eta) \;(g(\eta)\equiv g(y))$ by the MPS. In order to determine the intersection points  $\eta^{j}=(\eta_1^{j},\eta_{2}^{j},\eta_{3}^{j}) (j=\overline {1,N})$ of Wiener process trajectory and the surface $S_{1}$, we operate the following way.

 During the implementation of Wiener process, for each current point $\xi(t_k)$, defined from (3.3), its location with respect to $S_{1}$ is checked, i.e., for the point $\xi(t_k)$ the value
 $$
d=\xi_{1}^2+\xi_{2}^2+\xi_{3}^2
$$
is calculated and the following conditions: $d=1, d<1, d>1$ are checked. If $d=1$ then $\xi(t_k)\in S_{1}$ and $\eta^{j}=\xi(t_k)$. If $d<1$ then $\xi(t_k)\in K_{1}\backslash S_{1}$ and the process continues until its trajectory crosses the sphere $S_{1}$, and if $d>1$ then $\xi(t_k) \overline \in K_{1}$.

Let $\xi(t_{k-1})\in K_{1}\backslash S_{1}$ for the moment $t=t_{k-1}$ and $\xi(t_k)\overline \in K_{1}$ for the moment $t=t_k$. For determination of the point $\eta^j$, a parametric
equation of a line $L$ passing through the points $\xi(t_{k-1})$ and
$\xi(t_k)$ is first obtained:
$$
\left\{\begin {array}{ll}
\xi_1=\xi_1^{k-1}+(\xi_1^k-\xi_1^{k-1})\theta,\\&\\
\xi_2=\xi_2^{k-1}+(\xi_2^k-\xi_2^{k-1})\theta,\\&\\
\xi_3=\xi_3^{k-1}+(\xi_3^k-\xi_3^{k-1})\theta,\\&\\
\end {array}\right.
\eqno (4.3)
$$
where $(\xi_1,\xi_2,\xi_3)$ is a point of $L$, $\theta$ is a
parameter $(-\infty< \theta< \infty)$, and $\xi_i^{k-1} \equiv
\xi_i(t_{k-1}), \xi_i^k \equiv \xi_i(t_k) \;(i=1,2,3)$. After this , for
definition of the intersection points $\eta^*$ and $\eta^{**}$ of line $L$
and the sphere $S_{1}$ is solved an equation with respect to $\theta$.

It is easy to see that for the parameter $\theta$ we obtain an equation
$$
A\theta^2+2B\theta+C=0, \eqno (4.4)
$$
whose discriminant $d^*=B^2-AC>0.$

Since the discriminant of (4.4) is positive, the equation (4.4) has two solutions $\theta_1$ and $\theta_2$. Respectively, the points  $\eta^*$ and
$\eta^{**}$ are defined from the (4.3). In the role of the
point $\eta^j$ we choose the one from $\eta^*$ and $\eta^{**}$ for which
$|\xi(t_k)-\eta|$ is minimal.

In both examples, considered by us for the determination the
intersection points $\eta^j=(\eta_1^j, \eta_2^j, \eta_3^j) (j=\overline{1,N})$
of the trajectory of the Wiener process and the surface $S_{1}$ is used
that scheme, which is above described. As noted above, for
verification first we solve the auxiliary Problem $B$ with
the program of Problem $A^{*}$. In the numerical experiments for the example 1,
in test Problem B, we took $\xi^{0}=(0,5,0).$

In Table 4.1B the absolute errors $\Delta^i$ of the approximate
solution $v_N(\xi)$ of the test problem $B$ at the points $\xi^i \in K_1 \backslash S_{1}
(i=\overline{1,5})$ are presented. The notation $(E \pm k)$ is used for
$10^{\pm k}$.

\begin  {center}
\begin {tabular}{|c|c|c|c|c|c|}
\multicolumn{6}{l}{Table \ 4.1B. Results for Problem B (in Example 1 )}
\vspace {0.1mm}\\
\hline
$\xi^i$&$(0,0,0)$&$(0,0,0.5)$&$(0,0,-0.5)$&$(0,0,0.8)$&$(0,0,-0.8)$\\
\hline
$N$&$\Delta^1$&$\Delta^2$&$\Delta^3$&$\Delta^4$&$\Delta^5$ \\
\hline
$5E+3$&$0.47E-3$&$0.40E-3$&$0.85E-5$&$0.28E-3$&$0.23E-3$ \\
$1E+4$&$0.28E-3$&$0.25E-3$&$0.10E-3$&$0.13E-3$&$0.41E-4$ \\
$5E+4$&$0.62E-4$&$0.17E-3$&$0.83E-4$&$0.41E-4$&$0.28E-4$ \\
$1E+5$&$0.57E-4$&$0.44E-4$&$0.45E-4$&$0.36E-4$&$0.33E-4$ \\
$5E+5$&$0.54E-4$&$0.56E-4$&$0.28E-4$&$0.20E-5$&$0.35E-4$ \\
$1E+6$&$0.29E-5$&$0.43E-4$&$0.30E-4$&$0.51E-4$&$0.28E-4$ \\
\hline
\end {tabular}
\end {center}
\

On the basis of results presented in the Table 4.1B we can conclude that
the program for Problem $A^{*}$ is correct.

We conducted the check experiment. Namely, we calculated the
probabilistic solution of Problem $B$ at the point (0,0,0) for
$N=1E+5, nq=200$ and we obtained $\Delta^1=0.18E-4$ (see Table
4.1B). The result is improved, as expected (see Section 3). In general, if
more accuracy is need, then calculations for sufficiently large
values of $nq$ and $N$ must be performed. For this, it is appropriate to use parallel computing systems.

\begin  {center}
\begin {tabular}{|c|c|c|c|c|c|}
\multicolumn{6}{l} {Table \ $4.1A^{*}$. Results for Problem $A^{*}$ (in Example 1)}
\vspace {0.5mm}\\
\hline
$\xi^i $&$(0,0,0)$&$(0,0,0.5)$&$(0,0,-0.5)$&$(0,0,0.8)$&$(0,0,-0.8)$\\
\hline
$N$&$v_N(\xi^1)$&$v_N(\xi^2)$&$v_N(\xi^3)$&$v_N(\xi^4)$&$v_N(\xi^5)$ \\
\hline
$5E+3$&$1.50160$&$1.70600$&$1.69820$&$1.89940$&$1.90000$ \\
$1E+4$&$1.49500$&$1.69480$&$1.69440$&$1.89710$&$1.89430$ \\
$5E+4$&$1.50448$&$1.70062$&$1.69948$&$1.89780$&$1.89870$ \\
$1E+5$&$1.49664$&$1.69822$&$1.70120$&$1.89639$&$1.89609$ \\
$5E+5$&$1.50069$&$1.69935$&$1.69917$&$1.89577$&$1.89609$ \\
$1E+6$&$1.50009$&$1.69903$&$1.69849$&$1.89517$&$1.89587$ \\
\hline
\end {tabular}
\end {center}
\
In Table $4.1A^{*}$ the values of the approximate solution $v_N(\xi)$ to the
Problem $A^{*}$ at the same points $\xi^i (i=\overline {1,5})$ are given.
The boundary function (4.2) is symmetric with respect to the plane $Ox_1x_2$.
Respectively, in the role of $\xi^2, \xi^3$
and $\xi^4, \xi^5$, the points, which are symmetric with respect to the plane
$Ox_1x_2$, are taken. The results have sufficient accuracy for many
practical problems and agrees with the real physical picture.

\begin {center}
\begin {tabular}{|c|c|c|c|}
\multicolumn{4}{l}{Table \ $4.1A$.
 Results for Problem $A$ (in Example 1)}
\vspace{0.1mm}\\
\hline
$x^i$&$\xi^i$&$v_{N}(\xi^{i})$&$u_{N}(x^{i})$\\
\hline
$\infty$&$(0,0,0)$&$1.50069$&$0$\\
$(0,5000,5000)$&$(0,0.0001,0.0001$&$1.50038$&$2.12E-4$\\
$(0,1000,1000)$&$(0,0.0005,0.0005)$&$1.49971$&$0.00106$\\
$(0,100,100)$&$(0,0.005,0.005)$&$1.49957$&$0.01060$\\
$(0,50,50)$&$(0,0.01,0.01)$&$1.49990$&$0.02121$\\
$(0,10,10)$&$(0,0.05,0.05)$&$1.50153$&$0.10617$\\
$(0,1,1)$&$(0,0.5,0.5)$&$1.67766$&$1.18628$\\
\hline
\end {tabular}
\end {center}
\

In Table 4.1A the values of approximate solution $u_N(x)$ to Problem $A$ at the points
$x^{i}\in D (i=\overline {1,7})$ are given. For definition the values $u_N(x^i)$ 
the formula (2.9) is applied, respectively, in the role of points $x^i$ and $\xi^i (i=\overline{1,7})$, the points, which are symmetric with respect to the sphere $S_1$, are taken. It is evident that the points $x^i$ are situated on the ray, whose beginning is at the point $\xi=(0,0,0)$ and located in the plane $Ox_2x_3$. The ray's inclination angle with plane $Ox_1x_2$ is equal to $\pi/4.$
The points $\xi^i (i=\overline {1,7})$ are obtained by the inversion (2.1). The values $v_N(\xi^i)$ are the approximate solution $v_N(\xi)$ to Problem $A^{*}$ at the points $\xi^i$ for $N=1E+5$ and $nq=100.$
The obtained results agrees with the real physical picture.

As noted above, the program for Problem $A^{*}$ (in Example 1) is correct, therefore, in following example we can solve to Problem $A^{*}$ directly.
\\
\\
{\bf{ Example 2.}} Here in a role of infinite domain $D$ we took again
the exterior of unit sphere $S_{1}: y_{1}^2+y_{2}^2+y_{3}^2=1$ with the center at origin $O(0,0,0)$, where
$y(y_{1},y_{2},y_{3})$ is a point of the surface $S_{1}$. In considered case the boundary function $g(y)=g(y_{1},y_{2},y_{3})$ has the following form

$$
g(y)=\left\{\begin {array}{ll}
1,\;\;y\in S^1=\{y\in S_1|\;y_1>0,\;y_2>0,\;y_3>0\},\\&\\
0,\;\;y\in S^2=\{y\in S_1|\;y_1<0,\;y_2>0,\;y_3>0\},\\&\\
1,\;\;y\in S^3=\{y\in S_1|\;y_1<0,\;y_2<0,\;y_3>0\},\\&\\
0,\;\;y\in S^4=\{y\in S_1|\;y_1>0,\;y_2<0,\;y_3>0\},\\&\\
2,\;\;y\in S^5=\{y\in S_1|\;y_1>0,\;y_2>0,\;y_3<0\},\\&\\
1,\;\;y\in S^6=\{y\in S_1|\;y_1<0,\;y_2>0,\;y_3<0\},\\&\\
2,\;\;y\in S^7=\{y\in S_1|\;y_1<0,\;y_2<0,\;y_3<0\},\\&\\
1,\;\;y\in S^8=\{y\in S_1|\;y_1>0,\;y_2<0,\;y_3<0\},\\&\\
0,\;\;\;y\in l_k (k=1,2,3).
\end {array}\right .
\eqno (4.5)
$$

It is evident that in this case the discontinuity curves $l_k (k=1,2,3)$ are the circles, obtained by intersection of the coordinate planes and the sphere $S_1$. Actually, the sphere $S_1$ is divided into equal parts $S^i (i=\overline {1,8})$ by curves $l_k$. Besides, in the considered case, $l_k,S^2,S^4$ are non-conductors.

Since in this case the problem domain is the same that in Example 1, therefore for the determination of intersection points $\eta^{j}(j=\overline {1,N})$ of the Wiener process trajectory and the sphere $S_1$, the same algorithm, described in Example 1, is applied.

\begin {center}
\begin {tabular}{|c||c||c||c||c||c|}
\multicolumn{6}{l} {Table \ $ 4.2A^{*}$. Results for the Problem $A^{*}$ (in Example 2)}
\vspace {0.5mm}\\
\hline
$\xi^i $&$(0,0,0)$&$(0,0,0.8)$&$(0,0,-0.8)$&$(0.5,0.5,0.5)$&$(-0.5,-0.5,0.5)$\\
\hline
$N$&$v_N(\xi^1)$&$v_N(\xi^2)$&$v_N(\xi^3)$&$v_N(\xi^4)$&$v_N(\xi^5)$ \\
\hline
$5E+3$&$1.00580$&$0.56320$&$1.46020$&$0.97100$&$0.97180$ \\
$1E+4$&$0.99970$&$0.54620$&$1.45060$&$0.96960$&$0.06820$ \\
$5E+4$&$0.99992$&$0.55262$&$1.44466$&$0.96930$&$0.96774$ \\
$1E+5$&$0.99766$&$0.54962$&$1.45056$&$0.96976$&$0.96990$ \\
$5E+5$&$1.00060$&$0.55173$&$1.44925$&$0.97012$&$0.96989$ \\
$1E+6$&$1.00002$&$0.55116$&$1.44936$&$0.97000$&$0.97020$ \\
\hline
\end {tabular}
\end {center}
\
The values of the numerical solution of Problem $A^{*}$ at
the points $\xi^i \in K_1\backslash S_1 (i=\overline {1,5})$ are given in Table $4.2A^{*}$.
Since the boundary function (4.5) is symmetric with respect to the
axis $Ox_3$, therefore, in the role of
$\xi^i (i=4,5)$, the points which are symmetric with respect to the
axis $Ox_3$ are taken. The obtained results have sufficient
accuracy for many practical problems and agrees with the real physical picture.

\begin  {center}
\begin {tabular}{|c|c|c|c|}
\multicolumn{4}{l}{Table \ $4.2A.$ Results for Problem $ A $ (in Example 2)}
\vspace {0.5mm}\\
\hline
$x^i $&$\xi^{i}$&$v_{N}(\xi^{i})$&$u_{N}(x^{i})$\\
\hline
$\infty$&$(0,0,0)$&$1.00060$&$0$ \\
$(0,0,10^{4})$&$(0,0,10^{-4})$&$0.99961$&$9.9961E-5$ \\
$(0,0,10^{3})$&$(0,0,10^{-3})$&$1.00176$&$0.00102$ \\
$(0,0,10^{2})$&$(0,0,10^{-2})$&$0.99147$&$0.00915$ \\
$(0,0,10)$&$(0,0,10^{-1})$&$0.923912$&$0.09239$ \\
$(0,0,5)$&$(0,0,0.2)$&$0.85483$&$0.17097$ \\
$(0,0,4)$&$(0,0,0.25)$&$0.81979$&$0.20495$ \\
$(0,0,2)$&$(0,0,0.5)$&$0.67068$&$0.33534$ \\
$(0,0,1.5)$&$(0,0,2/3)$&$0.59743$&$0.39829$ \\
$(0,0,1.2)$&$(0,0,5/6)$&$0.54302$&$0.45252$ \\
\hline
\end {tabular}
\end {center}
\
The values of the numerical solution of Problem $A$ at the points $x^i \in D (i=\overline{1,10})$ are given in Table $4.2A.$ The values  $u_N(x^i)$ are obtained by the formula (2.9), therefore, in the role of points $x^i$ and $\xi^i (i=\overline {1,10})$, the points, symmetric to the sphere $S_1$ are taken. The points $x^i$ and $\xi^i$ are situated on the axis $Ox_3$, and the points $\xi^i$ are obtained by (2.1). The values $v_N(\xi^i)$ are numerical solution to Problem $A^{*}$
at the points $\xi^i (i=\overline {1,10})$ for $N=1E+5$ and $nq=100$. The obtained results are agrees with the real physical picture.

In this work we solved problems of type $A$  when the boundary functions $g_i(y) (i=\overline {1,m})$ are constants. This was motivated by our interest to find out how well the obtained results agree with the real physical picture. It is evident that solving Problem $A$ under condition (1.5) is as easy as the Problem $B$.

The analysis of results of numerical experiments show that the results
obtained by the proposed algorithm are reliable and it is effective
for the numerical solution of problems of type $A$ and $B$. In addition,
the algorithm is sufficiently simple for numerical implementation.
\\
\\
{\bf{6. Concluding remarks }}
\\
\\

1. In this work, we have demonstrated that the suggested algorithm is ideally suited for numerically solving of considered problems $A$ and $B$. It should be noted that it is possible using this algorithm to find the solution of a problem at any point of the domain, unlike other algorithms known in the literature.
\\

2. The algorithm does not require approximation of the boundary function,
which is one of its important properties.
\\

3. The algorithm is easy to program, its computational cost is low, it is
characterized by an accuracy which is sufficient for many practical problems.
\\

4. In the future we plan to investigate the following:
\\

  * Application of the proposed algorithm to the numerical solution of
  Dirichlet classical and generalized harmonic problems for the infinite space $R^3$ with a finite number
  of spherical cavities.
\\

  * Application of the MPS for the same type of problem for iregular pyramidal domains.
\\

  * Application of the MPS for the same type of problem in finite domains which are
  bounded by several closed surfaces.
\\

  * Application of the MPS for the same problem for infinite 2D domains with a finite
  number of circular holes.
\\

 \centerline {REFERENCES}
\vskip +0.3cm

1. G. A. Grinberg, The Selected Questions of Mathematical Theory of Electric
and  Magnetic Phenomena, Izd. Akad. Nauk SSSR, 1948,(in Russian).

2. W. R. Smythe, Static and Dynamic Electricity (second edition), New York, Toronto, London, 1950.

3. H. S. Carslaw, J. C. Jaegger, Conduction of Heat in Solids, Oxford University Press, London, 1959.

4. N. S. Koshlyakov, E. B. Gliner, M. M. Smirnov, Equations in Partial Derivatives
of Mathematical Physics, Moscow, 1970, (in Russian).

5. B. M. Budak, A. A. Samarski and A. N. Tikhonov, A Collection of Problems in
Mathematical physics, Nauka, Moscow, 1980, (in Russian).

6. L. V. Kantorovich and V. I. Krylov, Approximate Methods of Higher Analysis,
5-th edition, Fizmatgiz, Moscow--Leningrad, 1962,(in Russian).

7. N. Koblishvili, Z. Tabagari and M. Zakradze, On Reduction of the Dirichlet Generalized Boundary Value Problem to an Ordinary Problem for Harmonic Function. Proc. A.Razmadze Math. Inst., 132 (2003), 93-106.

8. M. Zakradze, N. Koblishvili, A. Karageorghis, Y. Smyrlis, On solving the dirichlet generalized problem for harmonic function by the method of fundamental solutions. Seminar of I. Vekua Institute of Applied Mathematics, REPORTS. 34 (2008), 24-32.

9. M. Kublashvili, Z. Sanikidze, M. Zakradze, A method of conformal mapping for solving the generalized dirichlet problem of Laplace's equcation. Proc. A.Razmadze Math. Inst., 160 (2012), 71-89.

10. A. Sh. Chaduneli, M. V. Zakradze and Z. A. Tabagari, A method of probabilistic
solution to the ordinary and generalized plane Dirichlet problem for the Laplace equation, Science and Computing, Proc. Sixth ISTC Scientific Advisory
Committee Seminar, Vol. 2, Moscow (20003), 361-366.

11. A. Chaduneli, Z. Tabagari, M. Zakradze, A computer simulation
of probabilistic solution to the Dirichlet plane boundary problem for the Laplace
equation in case of an infinite plane with a hole. Bull. Georg. Acad.Sci. Vol.171. No 3, 2005, 437-440.

12. A. Bayliss, M. Gunzburger, E. Turkel, Boundary Conditions for the Numerical Solution of Elliptic Equations in Exterior Regions. SIAM Journal on Applied Mathematics, 42(2), (1982), 430–451. 

13. Belgacem, Faker Ben, et al., On the Schwarz algorithms for the elliptic exterior boundary value problems. ESAIM: Mathematical Modelling and Numerical Analysis - Modélisation Mathématique et Analyse Numérique 39.4 (2005), 693-714. 

14. E. A. Kanunnikova, Analytic-numerical method to solve 3D exterior boundary problems for elliptic equations. St. Petersburg Polytechnic University Journal - Physics and Mathematics, 1 (2014), 35-38.

15. Antonios Charalambopoulos, Leonidas N. Gergidis, A novel stochastic method for the solution of direct and inverse exterior elliptic problems. Quarterly of Applied Mathematics, 76 (10) (2017), Iss. 1, DOI 10.1090/qam/1480.

16. J. A. Stratton, Electromagnetic Theory, McGraw-Hill, New York, 1941.

17. H. Lamb, Hydrodynamics, Springer-verlag,1967.

18. O. D. Kellogg, Foundations of Potential Theory, Dover, New York, 1953.

19. A. N. Tixonov, A. A. Samarskii, The Equations of Mathematical Physics, Nauka, Moscow, 1972,(in Russian).

20. M. Zakradze, M. Kublashvili, Z. Sanikidze, N. Koblishvili, Investigation and
numerical solution of some 3D internal Dirichlet generalized harmonic problems in finite domains, Transactions of A. Razmadze Mathematical Institute 171(2017)103-110.

21. E. B. Dynkin, A. A. Yushkevich, Theorems and problems on Markov's processes.
Nauka, Moscow, 1967 (in Russian).

22. A. Chaduneli, Z. Tabagari, M. Zakradze, On Solving the internal three-dimensional Dirichlet problem for a harmonic function by the method of probabilistic solution. Bull. Georg. Natl. Acad.Sci, (N.S) 2(1)(2008), 25-28.

23. M. Zakradze, Z. Sanikidze, Z. Tabagari, On solving the external three-dimensional Dirichlet problem for a harmonic function by the probabilistic method. Bull. Georg. Natl. Acad. Sci. (N.S) 4(3) (2010),19-23.

24. M. Zakradze, B. Mamporia, M. Kublashvili, N. Koblishvili , The method of probabilistic solution for 3D Dirichlet ordinary and generalized harmonic problems in finite domains bounded with one surface, Transactions of A. Razmadze Math. Inst. 172(2018),453-465.

25. M. Zakradze, M. Kublashvili, N. Koblishvili, A. Chakhvadze, The
method of probabilistic solution for determination of electric and
thermal stationary fields in conic and prismatic domains, Transactions of A. Razmadze Math. Inst. 174 (2020), issue 2, 235-246.

26. M. Zakradze, M. Kublashvili, Z. Tabagari, N. Koblishvili,
On numerical solving of the Dirichlet generalized harmonic problem for
regular $n$-sided pyramidal domains by the probabilistic method,
Transactions of A. Razmadze Math. Inst. 176 (2022), issue 1, 123-132.

27. M. Zakradze, Z. Tabagari, Z. Sanikidze, E. Abramidze, Computer modeling of a probabilistic solution for the Dirichlet generalized harmonic problem in some finite axisymmetric bodies with a cylindrical hole Bulletin of TICMI(the paper is in print).

28. V. S. Vladimirov, The Equations of Mathematical Physics. Nauka, Moscow, 1971 (in Russian).

\end {document}